\documentclass[11pt]{amsart}

\usepackage{graphicx}
\usepackage{amsmath, amsfonts, amssymb, amsthm}

\usepackage{color, hyperref}
\usepackage{fullpage}
\newtheorem{thm}{Theorem}[section]
\newtheorem{prop}[thm]{Proposition}
\newtheorem{lem}[thm]{Lemma}

\newtheorem{conj}[thm]{Conjecture} 
\theoremstyle{definition}

\theoremstyle{remark}
\newtheorem{remark}[thm]{Remark}

\newcommand{\qan}[2]{({#1}; q)_{#2}}
\newcommand{\qa}[1]{({#1}; q)_{\infty}}
\newcommand{\qna}[2]{({#1}; q^{#2})_{\infty}}

\begin{document}

\title[Proof of a Limited Version of Mao's Partition Rank Inequality]{Proof of a Limited Version of Mao's Partition Rank Inequality using a Theta Function Identity}

\author{Rupam Barman}
\address{Department of Mathematics, Indian Institute of Technology, Hauz Khas, New Delhi 110016, India}
\email{rupam@maths.iitd.ac.in}

\author{Archit Pal Singh Sachdeva}
\address{Department of Mathematics, Indian Institute of Technology, Hauz Khas, New Delhi 110016, India}
\email{antiarchit@gmail.com}

\subjclass[2010]{Primary: 11P83.}
\date{18th May, 2016}
\keywords{Partitions, ranks, rank differences, theta functions.}
\begin{abstract}
Ramanujan's congruence $p(5k+4) \equiv 0 \pmod 5$ led Dyson \cite{dyson} to conjecture the existence of a measure ``rank'' such that $p(5k+4)$ partitions of $5k+4$ could be divided into sub-classes with equal cardinality to give a direct proof of Ramanujan's congruence. The notion of rank was extended to rank differences by Atkin and Swinnerton-Dyer \cite{atkin}, who proved Dyson's conjecture. More recently, Mao proved several equalities and inequalities, leaving some as conjectures, for rank differences for partitions modulo 10 \cite{mao10} and for $M_2$ rank differences for partitions with no repeated odd parts modulo $6$ and $10$ \cite{maom2}. Alwaise et. al. proved four of Mao's conjectured inequalities \cite{swisher}, while leaving three open. Here, we prove a limited version of one of the inequalities conjectured by Mao. 
\end{abstract}

\maketitle

\section{Introduction and Results}

A \emph{partition} of a positive integer $n$ is a way of writing $n$ as a sum of positive integers, usually written in non-increasing order of the summands or parts of the partition. The number of partitions of $n$ is denoted by $p(n)$. For a partition $\lambda$, we denote the number of parts in the partition as $n(\lambda)$ and the largest part as $l(\lambda)$. 
\par
The celebrated Ramanujan's congruences for the partition function begged for a combinatorial interpretation: 
\begin{align*}
p(5k+4) &\equiv 0 \pmod{5}, \\
p(7k+5) &\equiv 0 \pmod{7}, \\
p(11k+6) &\equiv 0 \pmod{11}.
\end{align*}
Dyson \cite{dyson} defined the rank of a partition $\lambda$ to be $l(\lambda) - n(\lambda)$ and conjectured that partitions for $5k+4$ and $7k+5$ can be divided into five and seven equal sub-classes respectively based on their rank. Specifically, he claimed that
\begin{align*}
N(s, 5, 5n+4) &= \frac{p(5n+4)}{5}, \\ 
N(t, 7, 7n+4) &= \frac{p(7n+6)}{7},
\end{align*}
where $N(s, m, n)$ denotes the number of partitions of $n$ with rank $s$ modulo $m$. Atkin and Swinnerton-Dyer \cite{atkin} proved Dyson's conjecture by finding the generating functions for the rank differences $N(s, m, mk+d) - N(s, m, mk+d)$ for $k = 5, 7$. They obtained several other interesting identities apart from Ramanujan's congruences. 
\par 
Lovejoy and Osburn \cite{m2rankdiff} expanded on the work by Atkin and Swinnerton-Dyer to find rank differences for overpartitions and $M_2$ rank differences for partitions without repeated odd parts, which is defined for such a partition $\lambda$ by $ \left \lceil \frac{l(\lambda)}{2} \right \rceil - n(\lambda) .$ The corresponding count for number of partitions of $n$ with no repeated odd parts having its $M_2$ rank congruent to $s$ modulo $m$ is given by $N_2(s, m, n)$. They obtained all the rank difference formulas corresponding to $m = 3, 5$. 
\par 
Continuing on their work, Mao \cite{mao10, maom2} extended the results for Dyson rank differences modulo $10$ and $M_2$ rank differences modulo $6$ and $10$. He obtained several interesting inequalities based on his results such as
\begin{align*}
N(1,10,5n+1)  & > N(5,10,5n+1), \\
N_2(0,6,3n+1) + N_2(1,6,3n+1) & > N_2(2,6,3n+1) + N_2(3,6,3n+1).
\end{align*} 
\par 
Mao also gave some conjectures in \cite{mao10, maom2} based on computational evidence, both for the Dyson rank and $M_2$ rank for partitions with unique odd parts.
\begin{conj}
Computation evidence suggests that
\begin{align}
\label{mao10 conj a}
N(0, 10, 5n) + N(1, 10, 5n) &> N(4, 10, 5n) + N(5, 10, 5n),\\
\label{mao10 conj b}
N(1, 10, 5n) + N(2, 10, 5n) &\ge N(3, 10, 5n) + N(4, 10, 5n), \\
\label{mao610 conj b}
N_2(0, 10, 5n) + N_2(1, 10, 5n) &> N_2(4, 10, 5n) + N_2(5, 10, 5n),\\
\label{mao610 conj c} 
N_2(0, 10, 5n+4) + N_2(1, 10, 5n+4) &> N_2(4, 10, 5n+4) + N_2(5, 10, 5n+4),\\
\label{mao610 conj d}
N_2(1, 10, 5n) + N_2(2, 10, 5n) &> N_2(3, 10, 5n) + N_2(4, 10, 5n),\\
\label{mao610 conj e}
N_2(1, 10, 5n+2) + N_2(2, 10, 5n+2) &> N_2(3, 10, 5n+2) + N_2(4, 10, 5n+2), \\
\label{mao610 conj a} 
N_2(0, 6, 3n+2) + N_2(1, 6, 3n+2) &> N_2(2, 6, 3n+2) + N_2(3, 6, 3n+2).
\end{align}
In \eqref{mao10 conj b}, \eqref{mao610 conj d}, and \eqref{mao610 conj e}, $n \ge 1$, whilst in the rest $n \ge 0$. 
\end{conj}
Alwaise et. al. \cite[Theorem 1.3]{swisher} proved four of these seven inequalities conjectured by Mao, namely \eqref{mao10 conj a}, \eqref{mao10 conj b}, \eqref{mao610 conj b}, and \eqref{mao610 conj c} by using elementary methods based on the number of solutions of Diophantine equations solving for the exponents in the generating functions in the corresponding rank differences. They also observed that in \eqref{mao10 conj b}, the strict inequality holds. However, their methods weren't strong enough to prove the remaining three conjectures, which are still open. Here, we prove a limited version of \eqref{mao610 conj a}.
\begin{thm}\label{main}
Mao's conjecture \eqref{mao610 conj a} is true when $3 \nmid n + 1$. Specifically, we have that the following inequalities are true for all $n \ge 0$:
\begin{align}
N_2(0, 6, 9n+2) + N_2(1, 6, 9n+2) &> N_2(2, 6, 9n+2) + N_2(3, 6, 9n+2), \\
N_2(0, 6, 9n+5) + N_2(1, 6, 9n+5) &> N_2(2, 6, 9n+5) + N_2(3, 6, 9n+5).
\end{align}
\end{thm}
\section{Preliminaries}
The standard $q$-series notation is employed which is defined as
\begin{align*}
\qan{a}{n} &:= \prod_{i=0}^{n-1}(1-aq^i), \\
\qa{a} &:= \prod_{i=0}^\infty(1-aq^i),
\end{align*}
where $n \in \mathbb{N}$ and $a \in \mathbb{C}$. The empty product $\qan{a}{0}$ is defined to be $1$. 
\par 
The following elementary identities are used in manipulation of $q$-series to prove equalities between expressions. For $a, b \in \mathbb{Z}$, $c \in \mathbb{C}$, and for $k \in \mathbb{N}$, we have
\begin{align}
\label{oddunique}
\qa{-q} \cdot \qna{q}{2} &= 1, \\
\label{squares}
\qna{q^a}{b} \qna{-q^a}{b} &= \qna{q^{2a}}{2b}, \\
\label{bisection}
\qna{cq^a}{2b}\qna{cq^{a+b}}{2b} &= \qna{cq^a}{b}, \\
\label{ksection}
\qna{cq^a}{kb}\cdots\qna{cq^{a+(k-1)b}}{kb} &= \qna{cq^a}{b}.
\end{align} 
Further, we make use of the shorthand notation as employed by both Mao \cite{mao10, maom2} and Alwaise et. al. \cite{swisher}.
\begin{align*}
\qan{a_1, \dots, a_k}{n} &:= \qan{a_1}{n} \cdots \qan{a_k}{n}, \\
\qa{a_1, \dots, a_k} &:= \qa{a_1} \cdots \qa{a_k}, \\
J_b &:= \qna{q^b}{b}, \\
J_{a, b} &:= \qna{q^a, q^{b-a}, q^b}{b}. \\
\end{align*}
We will also use Mao's $M_2$ rank difference generating function to prove our result Theorem \ref{main}. Mao proved the following theorem which encapsulates the pertinent rank differences. 
\begin{thm}[Mao \cite{maom2}]\label{maogf} We have
\begin{align*}
&\sum_{n \ge 0} \left ( N_2(0, 6, n) + N_2(1, 6, n) - N_2(2, 6, n) - N_2(3, 6, n) \right )q^n \\
&= \frac{1}{J_{9, 36}} \sum_{n = - \infty}^{\infty} \frac{(-1)^nq^{18n^2+9n}}{1-q^{18n+3}} + q \frac{J_{6, 36}^2J_{18, 36}J_{36}^3}{J_{3, 36}^2J_{9, 36}J_{15, 36}^2}  \\
&+  \frac{J_{6, 36}J_{18, 36}^2J_{36}^3}{2qJ_{3, 36}^2J_{9, 36}J_{15, 36}^2} - \frac{1}{J_{9, 36}} \sum_{n=-\infty}^{\infty} \frac{(-1)^nq^{18n^2+9n-1}}{1+q^{18n}}.
\end{align*}
\end{thm}
Apart from this, an identity of Ramanujan theta function is also used. The Ramanujan's general theta function $f(a, b)$ is defined as 
\begin{align*}
f(a, b) &:= \sum_{n = -\infty}^{\infty} a^{\frac{n(n+1)}{2}}b^{\frac{n(n-1)}{2}} 
= (-a, -b, ab; ab)_\infty
\end{align*} 
with $|ab| < 1$ where the equality following through (and being equivalent to) Jacobi triple product identity. We will use the following two special cases of the theta function and the function $\chi(q)$ which are defined as
\begin{align}
\label{thetaid}
\varphi(q) &:= f(q, q) = \qna{-q, -q, q^2}{2}, \\
\label{psiid}
\psi(q) &:= f(q, q^3) = \frac{\qna{q^2}{2}}{\qna{q}{2}}, \\
\chi(q) &:= \qna{-q}{2}.
\end{align} 
The following theta function identity is used in the proof of our main result.
\begin{thm}[Baruah and Barman \cite{barman}]\label{theta}
We have \[ \varphi^2(q) + \varphi^2(q^3) = 2\varphi^2(-q^6)\frac{\chi(q)\psi(-q^3)}{\chi(-q)\psi(q^3)}. \]
\end{thm}
\section{Proof of Theorem \ref{main}}
We denote $d(n) := N_2(0, 6, n) + N_2(1, 6, n) - N_2(2, 6, n) - N_2(3, 6, n)$ for simplicity. We will show that the generating function $\sum_{n\ge0}d(3n+2)q^n$ has strictly positive coefficients for all $n \not\equiv 2 \pmod{3}$. We first compute the generating function $\sum_{n\ge0}d(3n+2)q^n$ using Theorem \ref{maogf}.
\begin{prop}\label{gfprop}
We have 
\begin{align*}\sum_{n \ge 0} d(3n+2) q^{n} = \frac{1}{qJ_{3,12}}  \left (\frac{J_{2, 12}J_{6, 12}^2J_{12}^3}{2J_{1, 12}^2J_{5, 12}^2} -  \sum_{n=-\infty}^{\infty} \frac{(-1)^nq^{6n^2+3n}}{1+q^{6n}} \right).
\end{align*}
\end{prop}
\begin{proof}
The proof is straightforward manipulation by including only exponents congruent to $2$ modulo $3$ in the original generating function, and then letting $q \mapsto
 q^{\frac{1}{3}}$ as follows:
 \begin{align*}
&\sum_{n \ge 0} d(3n+2) q^{3n+2}  = \frac{J_{6, 36}J_{18, 36}^2J_{36}^3}{2qJ_{3, 36}^2J_{9, 36}J_{15, 36}^2} - \frac{1}{J_{9, 36}} \sum_{n=-\infty}^{\infty} \frac{(-1)^nq^{18n^2+9n-1}}{1+q^{18n}} \\
\implies &\sum_{n \ge 0} d(3n+2) q^{3n} = \frac{J_{6, 36}J_{18, 36}^2J_{36}^3}{2q^3J_{3, 36}^2J_{9, 36}J_{15, 36}^2} - \frac{1}{J_{9, 36}} \sum_{n=-\infty}^{\infty} \frac{(-1)^nq^{18n^2+9n-3}}{1+q^{18n}} \\
\implies  &\sum_{n \ge 0} d(3n+2) q^{n} = \frac{J_{2, 12}J_{6, 12}^2J_{12}^3}{2qJ_{1, 12}^2J_{3, 12}J_{5, 12}^2} - \frac{1}{J_{3, 12}} \sum_{n=-\infty}^{\infty} \frac{(-1)^nq^{6n^2+3n-1}}{1+q^{6n}}.
\end{align*}
\end{proof}
\begin{remark}
Note that the while there is a $q$ in the denominator of the common factor above, it is canceled because the constant term of the expression inside the parentheses is zero. 
\end{remark}
We will also need the following lemma which will tie together the proof:
\begin{lem}\label{phiid}
We have \[ \frac{J_{2, 12}J_{6, 12}^2J_{12}^3}{J_{1, 12}^2J_{5, 12}^2} = \frac{\varphi^2(q) + \varphi^2(q^3)}{2}. \]
\end{lem}
\begin{proof}
We first write the expression in its constituent $q$-series and then use \eqref{squares} to cancel common factors in both numerator and denominator. We find that 
\begin{align*}
\frac{J_{2, 12}J_{6, 12}^2J_{12}^3}{J_{1, 12}^2J_{5, 12}^2} &= \frac{\qna{q^2, q^{10}, q^{12}}{12}\qna{q^6, q^6, q^{12}}{12}^2\qna{q^{12}}{12}^3}{\qna{q,q^{11}, q^{12}}{12}^2\qna{q^5,q^7, q^{12}}{12}^2} \\
&= \frac{\qna{q^2, q^{10}}{12}\qna{q^6, q^6, q^{12}}{12}^2}{\qna{q,q^{7}}{12}^2\qna{q^5,q^{11}}{12}^2} \\
&= \varphi^2(-q^6)\frac{\qna{q, q^{5}}{6}\qna{-q, -q^{5}}{6}}{\qna{q}{6}^2\qna{q^5}{6}^2} \\
&= \varphi^2(-q^6)\frac{\qna{-q, -q^5}{6}}{\qna{q, q^5}{6}}. \end{align*}
We next use \eqref{ksection} to reduce the $q$-series by multiplying the missing factors in both numerator and denominator, and simplify the expression based on \eqref{psiid} which is based on \eqref{oddunique}, to finally recognize the identity in Theorem \ref{theta} as follows:
\begin{align*}
\frac{J_{2, 12}J_{6, 12}^2J_{12}^3}{J_{1, 12}^2J_{5, 12}^2} &= \varphi^2(-q^6)\frac{\qna{-q, -q^5}{6}}{\qna{q, q^5}{6}} \\
&= \varphi^2(-q^6)\frac{\qna{-q, -q^5}{6}\qna{-q^3}{6}\qna{q^3}{6}}{\qna{q, q^5}{6}\qna{q^3}{6}\qna{-q^3}{6}} \\
&= \varphi^2(-q^6)\frac{\qna{-q}{2}\qna{q^6}{6}\qna{q^3}{6}}{\qna{q}{2}\qna{-q^3}{6}\qna{q^6}{6}} \\
&= \varphi^2(-q^6)\frac{\chi(q)\psi(-q^3)}{\chi(-q)\psi(q^3)} \\
&= \frac{\varphi^2(q) + \varphi^2(q^3)}{2} .
\end{align*}
\end{proof}
We now prove our result Theorem \ref{main}.
\begin{proof}[Proof of Theorem \ref{main}] We use Lemma \ref{phiid} and note that all the exponents of the infinite summation inside the parentheses are $0 \pmod 3$. Hence,
\begin{align*}
\sum_{n \ge 0} d(3n+2) q^{n} &= \frac{1}{qJ_{3,12}}  \left (\frac{J_{2, 12}J_{6, 12}^2J_{12}^3}{2J_{1, 12}^2J_{5, 12}^2} -  \sum_{n=-\infty}^{\infty} \frac{(-1)^nq^{6n^2+3n}}{1+q^{6n}} \right ) \\
&= \frac{1}{qJ_{3,12}} \left ( \frac{\varphi^2(q) + \varphi^2(q^3)}{4} - \frac{1}{2} + \sum_{n \ge 1} a_{3n}q^{3n} \right  ),
\end{align*}
where $a_{3n} \in \mathbb{Z}$.

Now let $3 \nmid n+1$, then
\begin{align*}
d(3n+2)  &= [q^{n}]  \frac{1}{qJ_{3,12}} \left ( \frac{\varphi^2(q) + \varphi^2(q^3)}{4} - \frac{1}{2} + \sum_{n \ge 1} a_{3n}q^{3n} \right  ) \\
&= [q^{n+1}] \left ( \frac{\varphi^2(q) + \varphi^2(q^3)}{4{J_{3,12}}} - \frac{1}{2} + \frac{1}{J_{3,12}} \sum_{n \ge 1} a_{3n}q^{3n} \right  ) \\
&= [q^{n+1}]  \frac{\varphi^2(q) + \varphi^2(q^3)}{4{J_{3,12}}} - [q^{n+1}] \frac{1}{2} + [q^{n+1}] \frac{1}{J_{3,12}} \sum_{n \ge 1} a_{3n}q^{3n} \\
&= [q^{n+1}]  \frac{\varphi^2(q) + \varphi^2(q^3)}{4{J_{3,12}}}
\end{align*}
where $[x^k]f(x)$ denotes the coefficient of $x^k$ in the generating function $f(x)$.
It now suffices to show that all coefficients of $\frac{\varphi^2(q) + \varphi^2(q^3)}{J_{3,12}}$ are positive. This follows as 
\begin{align*}
\frac{\varphi^2(q) + \varphi^2(q^3)}{J_{3,12}} &= \frac{2 + 4q + 4q^2 + \sum_{n\ge3}b_nq^n}{(1-q^3)\qna{q^9, q^9, q^{12}}{12}} \\
&= \left(2 + 4q + 4q^2 + \sum_{n\ge3}b_nq^n\right)\left(\sum_{n \ge 0}q^{3n}\right)\left(1 + \sum_{n\ge0}c_nq^n\right) 
\end{align*}
where $b_i$ and $c_i$ are non-negative. We can generate $q^{3n+k}$ using the above factors by $q^k$ from first, $q^{3n}$ from second, and $1$ from the last, where $k = 0, 1, 2$. This completes our proof for Theorem \ref{main}
\end{proof}

\section{Conclusion and Remarks}

The method employed by Alwaise et. al. doesn't work for this inequality because the expression inside the parentheses in Proposition \ref{gfprop} does seem to have negative coefficients for an infinite number of coefficients. 

This result is limited to $3n+2$ when $3 \nmid n+1$, but computational evidence suggests that $\displaystyle \frac{1}{1-q^{12}}  \left (\frac{J_{2, 12}J_{6, 12}^2J_{12}^3}{2J_{1, 12}^2J_{5, 12}^2} -  \sum_{n=-\infty}^{\infty} \frac{(-1)^nq^{6n^2+3n}}{1+q^{6n}} \right )$ has non-negative coefficients, and given the simplification with help of Lemma \ref{phiid}, a stronger version of the method used in \cite{swisher} along with using properties of $\varphi^2(q)$, in which coefficient of $q^n$ counts number of Diophantine solutions to $a^2 + b^2 = n$ might aid in proving the inequality when $3 \mid n + 1$.

\end{document}